\documentclass[11pt]{amsart}

\usepackage{amsmath,amsthm, amscd, amssymb, amsfonts}

\newcommand\boxee{\begin{tabular}{|p{0,3cm}|}
\hline \\ \hline \end{tabular}}

\newcommand\boxe{\begin{tabular}{|p{0,1cm}|}
\hline \\ \hline \end{tabular}}

\newcommand{\Vc}{{\mathcal V}}
\newcommand{\Hc}{{\mathcal H}}
\newcommand{\B}{{\mathcal B}}

\newcommand{\G}{{\mathcal G}}
\newcommand{\D}{{\mathcal D}}
\newcommand{\Ec}{{\mathcal E}}
\newcommand{\Pc}{{\mathcal P}}
\newcommand{\fde}{{\triangleright}}
\newcommand{\gde}{{\triangleleft}}

\newcommand{\Tc}{ {\mathcal T}}

\newcommand{\Ss}{{\mathcal S}}

\newcommand{\Bg}{{\mathcal B}}

\newcommand{\ot}{{\otimes}}
\newcommand{\sou}{\mathfrak s}
\newcommand{\tgt}{\mathfrak e}

\newcommand{\Fun}{\mbox{\rm Fun\,}}

\newcommand{\Do}{{\mathcal D}}
\newcommand{\uno}{{\bf 1}}
\newcommand{\ca}{{\mathcal C}}

\newcommand{\otb}{{\bar{\otimes}}}
\newcommand{\ku}{{\Bbbk}}

\newcommand{\id}{\mbox{\rm id\,}}

\newcommand{\op}{\mbox{\rm op\,}}

\newcommand{\Aut}{\mbox{\rm Aut\,}}

\newcommand{\Vect}{\mbox{\rm Vect\,}}

\newcommand\Rep{\operatorname{Rep}}
\newcommand\Opext{\operatorname{Opext}}

\theoremstyle{plain}
\numberwithin{equation}{section}

\newtheorem{teo}{Theorem}[section]
\newtheorem{lema}[teo]{Lemma}
\newtheorem{cor}[teo]{Corollary}
\newtheorem{prop}[teo]{Proposition}
\newtheorem{claim}{Claim}[section]

\theoremstyle{definition}
\newtheorem{defi}[teo]{Definition}
  \newtheorem{exa}[teo]{Example}
\theoremstyle{remark}

\newtheorem{rmk}[equation]{Remark}
\def\pf{\begin{proof}}
\def\epf{\end{proof}}
\theoremstyle{remark}

\begin{document}

\title[Tensor categories and Vacant Double Groupoids]{Tensor categories and Vacant Double Groupoids}
\author[Mombelli and Natale]{Juan Mart\'\i n Mombelli and Sonia Natale}
\thanks{This work was partially supported by
Agencia C\'ordoba Ciencia, ANPCyT-Foncyt, CONICET, Fundaci\'on
Antorchas and Secyt (UNC)}
\address{Facultad de Matem\'atica, Astronom\'\i a y F\'\i sica
\newline \indent
Universidad Nacional de C\'ordoba
\newline
\indent CIEM -- CONICET
\newline \indent Medina Allende s/n
\newline
\indent (5000) Ciudad Universitaria, C\'ordoba, Argentina} \email{
mombelli@mate.uncor.edu} \email{natale@mate.uncor.edu,
\newline \indent \emph{URL:}\/ http://www.mate.uncor.edu/natale}

\begin{abstract} We show that
fusion categories $\Rep(\ku^{\sigma}_{\tau} \Tc)$ of
representations of the weak Hopf algebra coming from a vacant
double groupoid $\Tc$ and a pair $(\sigma, \tau)$ of compatible
2-cocyles are group-theoretical.
\end{abstract}
\date{September 2, 2005}
\maketitle

\section{Introduction}

Semisimple tensor categories appear encoding symmetries of
distinct mathematical structures. This makes the problem of their
classification both a highly interesting and difficult one. The
problem is difficult even for special classes of these categories
like representations of finite dimensional semisimple Hopf
algebras. The main goal of this paper is to establish a relation
between two recent constructions of semisimple tensor categories.

\medbreak An important class of fusion categories was introduced
by Ostrik in \cite{O2}, also studied in the paper \cite{ENO}, by
Etingof, Nikshych and Ostrik. These fusion categories are built up
from finite group data, and are called {\it group-theoretical}.

Let $D$ be a finite group and let $V \subseteq D$ be a subgroup.
Let also $\omega\in Z^3(D,\ku^{\times})$ be a normalized
3-cocycle, $\psi\in C^2(V,\ku^{\times})$ a normalized 2-cochain,
such that $\omega|_{V\times V\times V} = d\psi$. Denote by $\ca(D,
\omega) = \Vect^D_{\omega}$ the tensor category of finite
dimensional $D$-graded vector spaces with associativity constraint
given by the 3-cocycle $\omega$. The twisted group algebra
$\ku_{\psi} V$ is an algebra in this category. The category
$\ca(D, \omega, V, \psi)$ is defined to be the category of
$\ku_{\psi} V$-bimodules in $\Vect^D_{\omega}$. This is a fusion
category with tensor product $\otimes_{\ku_{\psi} V}$ and unit
object $\ku_{\psi} V$. A fusion category is called
group-theoretical if it is equivalent to a category of the form
$\ca(D, \omega, V, \psi)$.

By the results in \cite{ENO}, the simple objects in a group
theoretical category have integer Frobenius-Perron dimensions and
thus every group theoretical category is the representation
category of a semisimple finite dimensional quasi-Hopf algebra. An
explicit description (up to gauge equivalence) of these quasi-Hopf
algebras was given in \cite{fs-indic}: the quasi-Hopf algebras
appearing in this description are a natural generalization of the
Dijkgraaf-Pasquier-Roche twisted quantum doubles $D^{\omega}G$
\cite{dpr}. It has been asked in \cite{ENO} if group-theoretical
Hopf algebras exhaust the class of semisimple Hopf algebras. The
answer to this question is not known until now.

\medbreak Another construction of a family of tensor categories
arising from finite vacant double groupoids was done in
\cite{AN1}, and later generalized to the nonvacant case in
\cite{AN2}. In these papers, a semisimple weak Hopf algebra
$\ku^{\vartheta}\Tc$ is naturally attached to a finite double
groupoid $\Tc$ satisfying a certain filling condition and a
certain perturbation datum $\vartheta$, and giving \textit{a
fortiori} a semisimple tensor category of representations.

\medbreak It is natural to ask if any fusion category arising from
this construction is group-theoretical. Several examples studied
in \cite{AN1} turn out to be group-theoretical. Necessary and
sufficient conditions on the double groupoid $\Tc$ in order that
the category $\Rep \ku \Tc$ be fusion  were given in \cite{AN2}.
It is also shown in \cite{AN2} that in the cases when $\vartheta$
comes from 'corner functions' intrinsically attached to the double
groupoid, the category $\Rep \ku^{\vartheta}\Tc = \Rep \ku\Tc$ is
the representation category of a finite-dimensional semisimple
quasi-Hopf algebra.

\medbreak In this paper we show that the tensor categories of the
form $\Rep \ku^{\sigma}_{\tau}\Tc$, arising from a \emph{vacant}
double groupoid $\Tc$ and a pair of compatible cocycles $\sigma,
\tau$ as in \cite{AN1}, are  group-theoretical whenever they are
fusion. See Theorems \ref{tequi} and \ref{gp-ttic}. The precise
condition under which these categories are fusion has been given
in \cite{AN1}: it reduces to the connectedness of the groupoid of
vertical edges. Our proof relies on the fact that vacant double
groupoids are essentially the same as \emph{matched pairs} of
groupoids \cite{mk1}. This allows us to prove generalizations of
certain category equivalences valid for finite groups, and get the
result.
 We give an explicit,
although not canonical, tensor equivalence between $\Rep
\ku^{\sigma}_{\tau}\Tc$ and a group-theoretical fusion category.
One of the main difficulties we encounter is that, although in the
connected case we may transfer 'group-theoretical facts' to
'groupoid-theoretical facts', some cohomological obstructions
appear. This is evidenced in Subsection \ref{connected}, c.f.
Proposition \ref{restriccion} therein.

It follows from our results that the Drinfeld center of the
category  $\Rep \ku^{\sigma}_{\tau}\Tc$ is equivalent to the
representation category of a twisted quantum double \cite{dpr}. In
view of the description in \cite{fs-indic} of group-theoretical
quasi-Hopf algebras, our results in the vacant context also bring
implicit a description of a quasi-Hopf algebra whose
representation category is equivalent to $\Rep
\ku^{\sigma}_{\tau}\Tc$.

\medbreak The paper is organized as follows: in Section
\ref{primera} we study several examples of tensor categories
arising from finite groupoids, which generalize group-theoretical
categories. In Section \ref{vacantes} we recall the definition and
characterizations of vacant double groupoids. Section
\ref{weak-hopf} contains the construction of the associated weak
Hopf algebras from \cite{AN1}; we also give in this section an
alternative description of the corresponding tensor category of
representations in terms of matched pairs of groupoids. Finally,
in Section \ref{main-res} we prove our main result.

\subsection*{Acknowledgement} The question we consider in this
paper arose from discussions with N. Andruskiewitsch. Both authors
thank him for encouraging.

\subsection*{Conventions and notation} We shall work over an
algebraically closed field $\ku$ of characteristic zero. The group
of units of $\ku$ is denoted by $\ku^{\times}$. All vector spaces
and algebras are assumed to be over $\ku$.

A groupoid with source $\sou$ and target $\tgt$ will be denoted as
a diagram $\tgt,\sou : \Do \rightrightarrows \Pc$. When $h, g \in
\Do$ are such that $\tgt(h) = \sou(g)$ their composition will be
indicated by $hg$. We shall identify the base $\Pc$ with a subset
of $\Do$ via the identity map: $P \mapsto \id_P \in \Do$.

We shall say that the groupoid $\Do \rightrightarrows \Pc$ is
\emph{finite} if $\Do$ (hence also $\Pc$) is a finite set. The set
of arrows $\alpha$ with source $\sou(\alpha) = P$ and target
$\tgt(\alpha) = Q$ will be indicated $\Do(P, Q)$; thus  $\Do(P,
Q)$ is a torsor over the group $\Do(P) = \Do(P, P)$. The groupoid
$\Do \rightrightarrows \Pc$ defines an equivalence relation $\sim$
on $\Pc$ by declaring $P \sim Q$ if and only if $\Do(P, Q) \neq
\emptyset$. We shall say that the groupoid is connected if this
relation is connected.

If $\Do \rightrightarrows \Pc$ is connected, then it is isomorphic
to the direct product $\Do(P) \times \Pc^2$, where $P$ is any
element of $\Pc$. In general, let $X$ be an equivalence class of
$\sim$ and let $\Do_X$ be corresponding connected groupoid  on the
base $X$ induced by $\Do$. Then $\Do \simeq \coprod_{X \in
\Pc/\sim} \Do_X$.

\medbreak

\section{Tensor categories coming from finite groupoids}\label{primera}

\subsection{Preliminaries on tensor categories}

The reader is referred to \cite{O1} for an exposition of the
notions of tensor categories and their module categories used
throughout in this paper.

\medbreak We shall be interested in a special class of tensor
categories, whose definition we recall next. These type of
categories have been intensively studied by a number of people in
the last years.

\begin{defi} A \emph{fusion} category over
$\ku$ is a $\ku$-linear semisimple and rigid tensor category $\ca$
such that:

\begin{itemize}\item[(i)] all hom spaces are finite dimensional;
\item[(ii)] the set of isomorphism classes of simple objects in
$\ca$ is finite; \item[(iii)] the unit object $\uno$ is simple.
\end{itemize}

\medbreak \noindent If $\ca$ satisfies (i) and (ii), then $\ca$ is
called a \emph{multi-fusion} category. See \cite{ENO}.
\end{defi}

\medbreak Suppose that $\ca$ is a multi-fusion category. Write
$\uno = \oplus_i\uno_i$, where $\uno_i$ are the simple
constituents of the unit object $\uno$. Then $\ca = \oplus_{i,
j}\ca_{ij}$, where $\ca_{ij}$ is the full subcategory whose
objects $X$ satisfy $X \otimes \uno_j = X = \uno_i \otimes X$. The
categories $\ca_{ii}$ are fusion categories called the
\emph{fusion components} of $\ca$.

Assume $\ca$ is indecomposable. Then $\mathcal M = \oplus_j
\ca_{ji}$ is an indecomposable module category over $\ca$
\cite{O1}, and there are isomorphisms \cite[2.4]{ENO}
$$\ca^*_{\mathcal M} = \Fun_{\ca}(\mathcal M, \mathcal M) \simeq
\ca_{ii}^{\op}.$$

\bigbreak \emph{For the rest of this section we fix a finite
groupoid $\tgt,\sou:\Do\rightrightarrows\Pc$ and let $\omega\in
Z^3(\Do, \ku^{\times})$ be a normalized 3-cocycle}; that is,
$\omega:\Do _{\tgt\!}\times_{\sou} \Do_{\tgt\!}\times_{\sou}\Do\to
\ku^{\times}$, such that
\begin{align*}\omega(a,b,c)\,\omega(a,bc,d)\,\omega(b,c,d)&=\omega(ab,c,d)
\,\omega(a,b,cd),\\
 \omega(a, P, b)&=1,\end{align*}
for all appropriately composable arrows $a,b,c,d \in \D$, $P \in
\Pc$.

\subsection{Pointed multi-fusion categories}\label{cat1}

Denote by $\ca(\Do,\omega)$ the category of finite dimensional
$\Do$-graded vector spaces with nontrivial associator $\omega$.
The category $\ca(\Do,\omega)$ is a multi-fusion category as
follows. If $M = \oplus_{\alpha\in\Do} M_{\alpha}$ and
$N=\oplus_{\beta\in\Do} N_{\beta}$ are two objects then $M\otimes
N := \oplus_{\tgt(\alpha) = \sou(\beta)}M_\alpha \otimes_{\ku}
N_\beta$ with $\Do$-grading determined by
$$(M\otimes N)_{\alpha}:=\bigoplus_{\alpha_1\alpha_2=\alpha}
M_{\alpha_1}\otimes_{\ku} N_{\alpha_2}.$$ The associativity
constraint on homogeneous elements is given by $\omega$. The unit
object is $ \ku\Pc$, with grading given by $|P|=\id_P$, $P\in
\Pc$. The associativity and unit isomorphisms for this category
are all canonical.

\medbreak The category $\ca(\Do,\omega)$ is a direct sum of its
indecomposable multi-fusion subcategories corresponding to
connected components of $\Do \rightrightarrows \Pc$. Therefore
$\ca(\Do,\omega)$ is not an \emph{indecomposable} multi-fusion
category unless $\Do \rightrightarrows \Pc$ is \emph{connected}.

\medbreak The irreducible constituents of $\uno = \ku\Pc$ are the
one-dimensional subspaces $\ku P$, $P \in \Pc$. Therefore,
$\ca(\Do,\omega)$ is a fusion category if and only if $\# \Pc =
1$, that is, if and only if $\Do$ is a \emph{group}. Indeed, the
categories $\ca_{P,Q}$ of $\ca = \ca(\Do,\omega)$ are the full
subcategories whose objects are finite dimensional vector spaces
graded by the set $\Do(P, Q)$ of arrows going from $P$ to $Q$.

\medbreak An object $X$ in a multi-fusion category $\ca$ will be
called \emph{invertible} if there exists an object $Y$ such that
$X \otimes Y$ and $Y \otimes X$ are irreducible summands of the
unit object $\uno$.

Following \cite{ENO} we shall say that a multi-fusion category is
\emph{pointed} if every simple object is invertible. Every pointed
multi-fusion category is equivalent to a category of the form
$\ca(\Do,\omega)$ for some finite groupoid $\Do$ and 3-cocycle
$\omega$.

\medbreak For all $P \in \Pc$, the component fusion category
$\ca_{P,P}$ is exactly the category of $\Do(P)$-graded vector
spaces with associativity given by the restriction
$\widehat\omega$ of $\omega$ to $\Do(P)$. In particular, every
fusion component in $\ca$ is \emph{pointed} group-theoretical.

Suppose that $\Do \rightrightarrows \Pc$ is connected. By
\cite{ENO}, $\mathcal M_{(P)} = \oplus_{Q \in \Pc}\ca_{QP}$ is an
indecomposable module category over $\ca$ consisting, in this
case, of all objects graded by arrows $\alpha$ with target
$\tgt(\alpha) = P$, such that $\ca^*_{\mathcal M_{(P)}} \simeq
\ca(\Do(P), \widehat \omega)^{\op}$.

\subsection{ The categories $\ca(\Do, \omega, \Vc, \psi)$}\label{cat2}
Let $\Vc \rightrightarrows \widetilde\Pc$ be a subgroupoid of
$\Do$ on the base $\widetilde\Pc \subseteq \Pc$. Let also $\psi:
\Vc {}_{\tgt}\times_{\sou}\Vc \to \ku^{\times}$ be a normalized
2-cochain such that
$$\omega \mid \Vc {}_{\tgt}\times_{\sou} \Vc
{}_{\tgt}\times_{\sou} \Vc = d\psi.$$ Then the twisted groupoid
algebra $\ku_{\psi}\Vc$  is an algebra in $\ca(\Do, \omega)$.

\begin{rmk} The category $\ca(\Do, \omega)_{\ku_{\psi}\Vc}$ of
$\ku_{\psi}\Vc$-modules in $\ca(\Do, \omega)$ is naturally a
module category over $\ca(\Do, \omega)$ \cite{O1}.

There is a decomposition $$\ca(\Do, \omega)_{\ku_{\psi}\Vc} =
\bigoplus_{P} {}_{(P)}\ca(\Do, \omega)_{\ku_{\psi}\Vc},$$ where
the module category ${}_{(P)}\ca(\Do, \omega)_{\ku_{\psi}\Vc}$ is
the full subcategory of $\ca(\Do, \omega)_{\ku_{\psi}\Vc}$ whose
objects are finite dimensional vector spaces graded by arrows with
target $P$. So that $\ca(\Do, \omega)_{\ku_{\psi}\Vc}$ is not
indecomposable.

In addition, the module category ${}_{(P)}\ca(\Do,
\omega)_{\ku_{\psi}\Vc}$ is indecomposable if and only if the
groupoid $\Vc \rightrightarrows \widetilde\Pc$ is
\emph{connected}.
\end{rmk}

\medbreak Denote by $\ca(\Do, \omega, \Vc, \psi)$ the tensor
category of $\ku_{\psi}\Vc$-bimodules in $\ca(\Do, \omega)$. The
category $\ca(\Do, \omega, \Vc, \psi)$ is thus the dual category
to $\ca (\Do, \omega)$ with respect to the module category
$\ca(\Do, \omega)_{\ku_{\psi}\Vc}$. By construction, the
categories $\ca(\Do, \omega, \Vc, \psi)$ are a multi-fusion
generalization of group-theoretical categories.

\begin{rmk}\label{equivalencias} Let $\eta: \Do {}_{\tgt}\times_{\sou} \Do \to
\ku^{\times}$, $\chi: \Vc \to \ku^{\times}$ be normalized
cochains. As in \cite[Remark 8.39]{ENO}, there are tensor
equivalences
$$\ca(\Do, \omega, \Vc, \psi) \simeq \ca(\Do, \omega', \Vc,
\psi'),$$ where $\omega' = \omega (d\eta)$, $\psi' = \psi \,
\eta\vert_{\Vc} \, (d\chi)$. \end{rmk}

Recall that a subgroupoid $\Vc \rightrightarrows \widetilde \Pc$
of $\Do$ is called \emph{wide} if $\widetilde \Pc = \Pc$.

\begin{lema} Suppose $\Vc$ is a wide subgroupoid. Then $\ca(\Do, \omega, \Vc, \psi)$ is a fusion category if and
only if $\Vc \rightrightarrows \Pc$ is connected. \end{lema}

\pf In this case, the connected components of $\Vc
\rightrightarrows \Pc$ correspond to irreducible constituents of
$\ku_{\psi}\Vc$ as an object in $\ca(\Do, \omega, \Vc, \psi)$.
This implies the lemma. \epf

Suppose that the 3-cocycle $\omega$ satisfies
$$\omega \mid \Vc {}_{\tgt}\times_{\sou} \Vc
{}_{\tgt}\times_{\sou} \Vc = 1.$$ Objects in $\ca(\Do,\omega,\Vc)$
can be described explicitly as follows: these are vector spaces
graded by arrows in $\Do$ with source and target in
$\widetilde\Pc$,
$$M = \oplus_{\alpha \in \Do\vert_{\widetilde\Pc}} M_{\alpha},$$
 together with actions of $\Vc$ by linear
isomorphisms:
$$g\rightharpoonup : M_\alpha \to M_{g\alpha}, \quad
\leftharpoonup h: M_\alpha \to M_{\alpha h}, $$ $g, h \in \Vc$,
$\alpha \in \Do$, such that $\tgt(g) = \sou(\alpha)$, $\sou(h) =
\tgt(\alpha)$.

\medbreak Letting $|m| \in \Do$ denote the homogeneity degree of
an homogeneous element $m \in M$, the following relations are
satisfied:
\begin{equation}\label{sw3} |m\leftharpoonup h|=|m|h, \,\,\,
|g\rightharpoonup m|=g|m|,
\end{equation}
\begin{align}\label{sw1}
\quad\quad gh\rightharpoonup m = \omega(g, h, |m|)
g\rightharpoonup( h\rightharpoonup m),
\end{align}
\begin{align}\label{sw4} \id_p\rightharpoonup m= \begin{cases}
m \quad &\text{if } p=\sou(|m|) \in \widetilde\Pc,\\
0 \quad &\text{otherwise, }
\end{cases}\\
\label{sw5} (m\leftharpoonup g)\leftharpoonup h=\omega(|m|,g,h)\,
m\leftharpoonup gh,\\
\label{sw6} m\leftharpoonup \id_p= \begin{cases} m \quad &\text{if } p=\tgt(|m|) \in \widetilde\Pc,\\
0 \quad &\text{otherwise, }
\end{cases}
\end{align}
\begin{equation}\label{sw2}
(g\rightharpoonup m)\leftharpoonup h = \omega(g, |m|, h)
g\rightharpoonup(m\leftharpoonup h),
\end{equation}
for any homogeneous element $m\in M$,  $g, h \in \Vc$, such that
the elements $g,|m|,h$ are composable in $\Do$.

\bigbreak

The tensor structure on the category $\ca(\Do,\omega,\Vc)$ can be
described as follows. If $M, N\in \ca(\Do,\omega,\Vc)$ denote by
$V(M,N)$ the (graded) subspace of $M\otimes N$ spanned by
\begin{equation}\label{tp1} (m\leftharpoonup g) \ot n- \omega(|m|,g,|n|) m \ot
(g\rightharpoonup n), \end{equation} where $m\in M$, $n\in N$, are
homogeneous elements such that $|m|,g,|n|$ are composable in $\D$.

Tensor product $\otb$ on $\ca(\Do,\omega,\Vc)$  is defined by
$M\otb N= (M\ot N)/ V(M,N)$, with $\Do$-grading inherited from $M
\otimes N$. The left and right actions of $\ku\Vc$ on $M\otb N$
are as follows. The class of $m\ot n$ will be denoted by $m\otb
n$. Assume that $m\in M, n\in N$ are homogeneous elements such
that $\tgt(|m|)=\sou(|n|)$. Let $g, h\in \Vc$ such that
$\tgt(g)=\sou(|m|), \tgt(|n|)=\sou(h),$ then
\begin{equation}\label{ap-asoc}\begin{aligned}g\rightharpoonup (m\otb n) & =  \omega^{-1}(g, |m|, |n|)
(g\rightharpoonup m)\otb n, \\
 (m\otb n)\leftharpoonup h & = \omega(|m|, |n|, h) m\otb (n\leftharpoonup h).\end{aligned}\end{equation}
The unit object is $\ku\Vc$.

\medbreak
\begin{exa} According to the above description, for all $P \in \Pc$, there is a natural
identification $$\ca(\Do, \omega, \Vc(P)) = \ca(\Do(P),
\widehat\omega, \Vc(P));$$ hence $\ca(\Do, \omega, \Vc(P))$ is a
fusion group-theoretical category. \end{exa}

\medbreak
\begin{exa} Let $P \in \Pc$, and consider the subgroupoid
$$\{ \id_P \} = \Vc_{(P)} \rightrightarrows \widetilde\Pc = \{ P
\}.$$ In this case, the category ${}_{\ku \Vc_{(P)}}\ca(\Do,
\omega)$ of $\ku\Vc_{(P)}$-modules in $\ca(\Do, \omega)$ is the
full tensor subcategory whose objects are graded by arrows
$\alpha$ with $\tgt(\alpha) = P$.

This coincides with the module category $\mathcal M_{(P)} =
\oplus_{Q \in \Pc} \ca_{QP}$ \cite[2.4]{ENO}.

\medbreak \emph{Suppose that the groupoid $\Do \rightrightarrows
\Pc$ is connected.} Let $D = \Do(P)$, and $\widehat\omega$ the
3-cocycle on $D$ obtained by restriction. As remarked before, by
the results in the proof of Proposition 2.17 in  \cite[5.5]{ENO},
$\mathcal M_{(P)}$ is an indecomposable module category and there
are isomorphisms
$$\ca^*_{\mathcal M_{(P)}} \simeq \ca(\Do, \omega, \Vc_{(P)}) \simeq
\ca(D, \widehat\omega)^{\op}.$$ The last isomorphism can also be
seen as a consequence of the description of $\ca(\Do, \omega,
\Vc_{(P)})$ given above.

In particular, the category $\ca(\Do, \omega, \Vc_{(P)})$ is a
fusion group-theoretical category, even though the groupoid
$\Vc_{(P)}$ may not be a \emph{wide} subgroupoid, that is, we may
have $\widetilde\Pc \neq \Pc$.
\end{exa}

\subsection{The connected case}\label{connected} From now on we shall
assume that $\Vc \rightrightarrows \widetilde\Pc$ is  a
\emph{wide} subgroupoid (so $\widetilde\Pc = \Pc$) which is
\emph{connected}. Thus $\Do \rightrightarrows \Pc$ is also
connected. We shall also fix a 3-cocycle $\omega$ on $\Do$ such
that $$\omega \mid \Vc {}_{\tgt}\times_{\sou} \Vc
{}_{\tgt}\times_{\sou} \Vc = 1.$$

The connectedness assumption on $\Vc$ implies that the unit object
$\ku\Vc$ in the category $\ca(\Do,\omega,\Vc)$ is irreducible.
Therefore, $\ca(\Do,\omega,\Vc)$ is in this case a \emph{fusion}
category.  In what follows we shall give an explicit equivalence
that shows that $\ca(\Do,\omega,\Vc)$ is indeed a
group-theoretical category.

\medbreak Fix an element $O\in\Pc$. Since $\Vc$ is connected, for
any $P\in\Pc$, we can choose an element $\tau_P\in \Vc(O,P)$.
There is no harm to assume $\tau_O = \id_O$.

\medbreak Set $V = \Vc(O), D = \D(O)$.  For all $P,Q \in \Pc$ we
have  $$\D(P,Q) = \tau^{-1}_P D \tau_Q, \quad \Vc(P,Q)=\tau^{-1}_P
V \tau_Q.$$ Let $\widehat{\omega}$ be the restriction of $\omega$
to $D$; so $\widehat{\omega}$ gives a well-defined 3-cocycle on
$D$. This restriction map determines a group isomorphism
$$H^3(\D,\ku^{\times})\to H^3(D,\ku^{\times}), \quad \omega\mapsto
\widehat{\omega}.$$ See for instance \cite{AM}.

\medbreak Let $\widetilde\omega \in Z^3(\Do, \ku^{\times})$ be the
3-cocycle defined by the formula
\begin{equation}\label{rest-omega}\widetilde\omega(g, h, t) = \widehat \omega(\tau_Pg\tau_Q^{-1}, \tau_Qh\tau_R^{-1},
\tau_Rt\tau_S^{-1}), \end{equation} for $g \in \Do(P, Q)$, $h \in
\Do(Q, R)$, $t \in \Do(R, S)$.

Note that, since $\tau_P \in \Vc$, for all $P \in \Pc$, then we
also have $\widetilde\omega\vert\Vc {}_{\tgt}\times_{\sou} \Vc
{}_{\tgt}\times_{\sou} \Vc = 1$.

\medbreak By construction, $\widetilde\omega\vert_{D} =
\widehat\omega$. Therefore $\widetilde\omega$ is cohomologous to
$\omega$, say \begin{equation}\label{coborde}\widetilde\omega =
\omega (d\psi),\end{equation} for some normalized 2-cochain $\psi:
\Do {}_{\tgt}\times_{\sou} \Do \to \ku^{\times}$.

In particular, the tensor categories $\ca(\Do, \omega)$ and
$\ca(\Do, \widetilde\omega)$ are equivalent.

\medbreak Moreover, since both restrictions $\omega\vert_{\Vc}$
and $\widetilde\omega\vert_{\Vc}$ are trivial, relation
\eqref{coborde} implies that the restriction $\psi: \Vc
{}_{\tgt}\times_{\sou} \Vc \to \ku^{\times}$ is a 2-cocycle, and
by Remark \ref{equivalencias} there is a tensor equivalence
\begin{equation}\label{prim-equiv}\ca(\Do, \omega, \Vc) \simeq
\ca(\Do, \widetilde\omega, \Vc, \psi).\end{equation}

Let $\widehat\psi = \psi\vert_V: V \times V \to \ku^{\times}$ be
the corresponding 2-cocycle on $V$, and let also $\widetilde\psi:
\Vc {}_{\tgt}\times_{\sou} \Vc \to \ku^{\times}$ be the 2-cocycle
defined by the formula
\begin{equation}\label{rest-psi}\widetilde\psi(g, h) : = \widehat\psi(\tau_Pg\tau_Q^{-1},
\tau_Qh\tau_R^{-1}), \end{equation} for all $g \in \Vc(P, Q)$, $h
\in \Vc(Q, R)$. Repeating the argument used above for $\omega$,
the results in \cite{AM} imply that $\psi$ and $\widetilde\psi$
differ by a 2-coboundary, that is,
$$\widetilde\psi = \psi (d\chi),$$ on the connected subgroupoid $\Vc$, where $\chi: \Vc \to
\ku^{\times}$ is a normalized 1-cochain.

Combining this fact with the equivalence \eqref{prim-equiv} and
Remark \ref{equivalencias} we get tensor equivalences
\begin{equation}\label{seg-equiv}\ca(\Do, \omega, \Vc) \simeq \ca(\Do, \widetilde\omega, \Vc, \psi) \simeq
\ca(\Do, \widetilde\omega, \Vc, \widetilde\psi). \end{equation}

\bigbreak Let $M \in \ca(\D, \widetilde\omega, \Vc,
\widetilde\psi)$. Let us define $F(M) \in \ca(D, \widehat{\omega},
V, \widehat\psi)$ by
$$F(M)=\bigoplus_{z\in D} M_z.$$
This gives us a well-defined functor $F: \ca(\D, \widetilde\omega,
\Vc, \widetilde\psi) \to \ca(D, \widehat{\omega}, V,
\widehat\psi)$.

Note that, for $M \in \ca(\D, \widetilde\omega, \Vc,
\widetilde\psi)$,
\begin{equation}\label{bimodule} M=\bigoplus_{P,Q\in\Pc}
\tau^{-1}_P\rightharpoonup F(M)\leftharpoonup\tau_Q.
\end{equation}

\begin{prop}\label{restriccion} The functor $F$ induces an equivalence of tensor categories
$$F: \ca(\D, \widetilde\omega, \Vc, \widetilde\psi) \overset{\simeq}\to
\ca(D, \widehat{\omega}, V, \widehat\psi).$$
\end{prop}

\begin{proof} Define a functor $G: \ca(D, \widehat{\omega}, V, \widehat\psi) \to \ca(\D,
 \widetilde\omega, \Vc, \widetilde\psi)$ by the formula
$$G(W)=\bigoplus_{P,Q\in\Pc}\bigoplus_{z\in D}
 W_{(P,z,Q)} \simeq \ku\Pc \otimes_{\ku}W \otimes_{\ku} \ku\Pc,$$
for all $P,Q\in\Pc, z\in D$, where $$W_{(P,z,Q)} =
\ku\tau^{-1}_P\ot_{\ku} W_z\ot_{\ku} \ku\tau_Q, $$ is an object in
the category $\ca(\D, \widetilde\omega,\Vc, \widetilde\psi)$ as
follows: all elements in $W_{(P,z,Q)}$ are homogeneous of degree
$\tau^{-1}_P z\tau_Q$. If $v\in W_z$ and $g\in\Vc(R,P),
h\in\Vc(Q,S)$ then the left and right actions are defined by
\begin{align*} g\rightharpoonup \tau^{-1}_P\ot w \ot\tau_Q&=
\tau^{-1}_R\ot (\tau_R g \tau^{-1}_P)\cdot w \ot\tau_Q \in W_{(R,z,Q)},\\
\tau^{-1}_P\ot w \ot\tau_Q\leftharpoonup h&= \tau^{-1}_P\ot w
\cdot (\tau_Q h \tau^{-1}_S)\ot\tau_S\in W_{(P,z,S)}.
\end{align*}

Formulas \eqref{rest-omega} and \eqref{rest-psi} imply that $G$ is
well defined.

\medbreak We claim that the functors $F$ and $G$ are inverse
tensor equivalences.

Because $\widetilde \omega\vert_{D} = \widehat \omega$ and
$\widetilde \psi\vert_{V} = \widehat \psi$, we have $F(G(W))
\simeq W$, for all $W \in \ca(D, \widehat\omega, V,
\widehat\psi)$.

The proof of the proposition will be finished once we have
established the following claim. Its proof uses connectedness of
the groupoid $\Vc \rightrightarrows \Pc$.

\begin{claim} The map $$f: (\ku\tau^{-1}_P\ot_{\ku}
U\ot_{\ku} \ku\tau_Q) \otimes_{\ku} (\ku\tau^{-1}_Q\ot_{\ku}
W\ot_{\ku} \ku\tau_S) \to \ku\tau^{-1}_P\ot_{\ku} \left(U
\overline\otimes W\right)\ot_{\ku} \ku\tau_S,$$ determined by
$$\tau^{-1}_P\ot_{\ku} u \ot_{\ku} \tau_Q \otimes_{\ku} \tau^{-1}_Q \ot_{\ku}
w \ot_{\ku}\tau_S \mapsto \tau^{-1}_P\ot_{\ku} (u \overline\otimes
w) \ot_{\ku}\tau_S,$$ induces a natural isomorphism $\zeta_{U,
W}:G(U)\otb G(W)\to G(U \otb W)$. This isomorphism endows $G$ with
a tensor functor structure.
\end{claim}

Recall that $\otimes$ means tensor product in the category of
$\Do$-graded spaces.

\pf[Proof of the Claim] First note that $f$ defines a map from the
space
$$(\ku\Pc \otimes_{\ku}U\otimes_{\ku}\ku\Pc) \otimes (\ku\Pc \otimes_{\ku}W\otimes_{\ku}\ku\Pc)
= \bigoplus_{P, Q, S}\ku\tau^{-1}_P\ot_{\ku} U\ot_{\ku} \ku\tau_Q
\otimes_{\ku} \ku\tau^{-1}_Q\ot_{\ku} W\ot_{\ku} \ku\tau_S$$ to
the space
$$G(U \overline\otimes W) = \bigoplus_{P, S}\ku\tau^{-1}_P\ot_{\ku} \left(U \overline\otimes W\right)\ot_{\ku}
\ku\tau_S.$$ Let now  $g \in \Vc$ with $\sou(g) = Q$, $\tgt(g) =
Q'$. Then we have
\begin{align*}& f\left((\tau^{-1}_P  \ot_{\ku}  u  \ot_{\ku} \tau_Q) \leftharpoonup g \otimes_{\ku} \tau^{-1}_{Q'}
\ot_{\ku} w \ot_{\ku} \tau_S \right)  \\
& = f(\tau^{-1}_P\ot_{\ku}
u.(\tau_Q g \tau_{Q'}^{-1}) \ot_{\ku} \tau_{Q'} \otimes_{\ku}
\tau^{-1}_{Q'} \ot_{\ku}
w \ot_{\ku}\tau_S)  \\
& = \tau^{-1}_P\ot_{\ku} \overline{u.(\tau_Q g \tau_{Q'}^{-1}) \ot_{\ku} w} \ot_{\ku}\tau_S \\
& = \widehat \omega (|u|, \tau_Q g \tau_{Q'}^{-1}, |w|) \,
\tau^{-1}_P \ot_{\ku} \overline{u \ot_{\ku} (\tau_Q g \tau_{Q'}^{-1}).w} \ot_{\ku} \tau_S  \\
& = \widetilde \omega(\tau_P|u|\tau_Q^{-1}, g,
\tau_{Q'}|w|\tau_S^{-1}) \, f\left(\tau^{-1}_P  \ot_{\ku} u
\ot_{\ku} \tau_Q \otimes_{\ku} g \rightharpoonup (\tau^{-1}_{Q'}
\ot_{\ku} w \ot_{\ku}\tau_S) \right).
\end{align*} Hence $f$ induces a map $\zeta_{U, W}:G(U)\otb G(W)\to G(U \otb
W)$. Using \eqref{ap-asoc} and \eqref{rest-omega} we see that
$\zeta_{U, W}$ is a map in $\ca(\D,
 \widetilde\omega, \Vc, \widetilde\psi)$.

Let $\overline f: \ku\Pc \ot_{\ku} U \ot_{\ku}W \ot_{\ku}\ku\Pc
\to G(U) \overline{\otimes} G(W)$ be determined by $$\overline
f(\tau^{-1}_P \ot_{\ku} u \ot_{\ku} w \ot_{\ku} \tau_S) =
\tau^{-1}_P \ot_{\ku} u \ot_{\ku} \tau_O \otimes \tau_O
\otimes_{\ku} w \ot_{\ku} \tau_S.$$ Recall that $\tau_O = \id_O$.
A similar computation to the one done before shows that $\overline
f$ induces a map $\xi_{U, W}: G(U \overline{\ot} W) \to G(U)
\overline{\otimes} G(W)$. Clearly $\zeta_{U, V}\xi_{U, V} = \id$.
On the other hand, for all $P, Q, S \in \Pc$, we have
\begin{align*}\tau^{-1}_P \ot_{\ku} u \ot_{\ku} \tau_Q  &
\overline{\otimes}
\tau_Q^{-1} \otimes_{\ku} w \ot_{\ku} \tau_S
\\ & = \left( \tau^{-1}_P \ot_{\ku} u \ot_{\ku} \tau_O \right)
\leftharpoonup \tau_Q \overline{\otimes} \tau_Q^{-1}
\rightharpoonup \left(
\tau_O \otimes_{\ku} w \ot_{\ku} \tau_S \right) \\
& = \tau^{-1}_P \ot_{\ku} u \ot_{\ku} \tau_O \overline{\otimes}
\tau_O \otimes_{\ku} w \ot_{\ku} \tau_S,\end{align*} because
$\widetilde \omega$ and $\widetilde \psi$ are trivial whenever any
one of their arguments is a $\tau_Q$. This implies that $\xi_{U,
V}\zeta_{U, V} = \id$. Therefore $\xi_{U, V}$ is inverse to
$\zeta_{U, W}$, whence $\zeta$ is an isomorphism.

The definition of $\widetilde\omega$ also implies that $\zeta$ is
compatible with the associativity and unit constraints. This
finishes the proof of the claim. \epf
\end{proof}

\begin{cor}\label{equiv-bimod} There is a tensor equivalence $\ca(\D, \omega, \Vc) \simeq \ca(D, \widehat \omega, V, \widehat \psi)$.
In particular, the category $\ca(\D, \omega, \Vc)$ is
group-theoretical. \end{cor}

\pf It follows from \eqref{seg-equiv} and Proposition
\ref{restriccion}. \epf

\section{Vacant double groupoids}\label{vacantes}

According to the definition given by Ehresmann \cite{ehr}, a {\it
double groupoid} is a groupoid object in the category of
groupoids. A double groupoid can be represented in the form of
four related groupoids
$$\Tc :\quad \begin{matrix} \B &\rightrightarrows &\Hc
\\\downdownarrows &&\downdownarrows \\ \Vc &\rightrightarrows &\Pc \end{matrix}$$
subject to a set of axioms.

This structure admits a pictorial description as 'boxes' that can
be composed in two directions: horizontal and vertical. Throughout
we shall work with the conventions and notations from
\cite[Section 2]{AN1}. An element $A\in \B$ is depicted as a box
$$A =
\begin{matrix} \quad t \quad \\ l \,\, \boxe \,\, r \\ \quad b\quad
\end{matrix},$$ where $t = t(A), b = b(A) \in \Hc$ are respectively the source and
target of $A$ with respect to vertical composition, and similarly
for $l = l(A), r = r(A) \in \Vc$ with respect to horizontal
composition.

\medbreak Horizontal and vertical composition of boxes will be
written from left to right and from top to bottom, respectively.
 The notation $AB$
(respectively, $\begin{matrix}A \vspace{-4pt}\\B\end{matrix}$)
will indicate the horizontal (respectively, vertical) composition;
this notation will always implicitly assume that $A$ and $B$ are
composable in the appropriate sense.

\medbreak A double groupoid $\Tc$ is called \emph{vacant} if for
any $g\in \Vc$, $x\in \Hc$ such that the target of $x$ coincides
with the source of $g$, there is exactly one  box $X \in \B$ such
that $X = \begin{matrix} \quad x \quad \\  \,\, \boxe \,\, g \\
\quad \quad \end{matrix}$.

In particular, in a vacant double groupoid, every box is
determined by (any) pair of adjacent edges.

\subsection{Matched pairs of groupoids}
Let $\sou, \tgt: \G \rightrightarrows \Pc$ be a groupoid with
source and target maps $\sou$ and $\tgt$. A {\it left action} of
$\G$ on a map $p:{\mathcal E}\rightarrow\Pc$ is a map $\fde:\G
{\,}_{\tgt}\times_p \Ec \to \Ec$ such that, for all composable
$g,h\in\G$, $x\in\Ec$,
\begin{equation*} p(g\fde x)=\sou(g),\qquad
 g \fde(h \fde x)=gh \fde x,\qquad\id_{p(x)} \, \fde  x =  x.
\end{equation*}

\medbreak Similarly, a {\it right action} of $\G$ on ${\mathcal
E}$ is a map $\gde :\Ec {\,}_p\times_{\sou} \G \to \Ec$ such that
\begin{equation*}p(x\gde g)=\tgt(g),\qquad
(x \fde g) \fde h=x \gde gh,\qquad x \gde \id_{p(x)}  =  x,
\end{equation*}
for all composable $g,h\in\G$, $x\in\Ec$.

\medbreak A {\it matched pair of groupoids} is a collection
$(\Hc,\Vc,\fde,\gde)$, where  $t, b:\Vc\rightrightarrows \Pc$ and
$l, r:{\mathcal H}\rightrightarrows \Pc$ are two groupoids on the
same base $\Pc$, $\fde:{\mathcal H}_{\,\tgt\!}\times_{\sou}
\Vc\rightarrow \Vc$ is a left action of $\Hc$ on $(\Vc, b)$ and
$\gde:{\mathcal H}_{\,\tgt\!}\times_{\sou} \Vc\rightarrow \Hc$ is
a right action of $\Vc$ on $(\Hc, l)$ satisfying the following
conditions:
\begin{equation*} \label{m7}  r(x\fde g) =
t(x\gde g),\quad x\fde gh =(x\fde g)((x\gde g)\fde h),\quad
 xy\gde g = (x\gde (y\fde g))(y\gde g),
\end{equation*}
for composable elements $x,y\in\Hc$ and $g,h\in\Vc$.

 \subsection{Exact factorizations} Let $(\Hc,\Vc, \gde,
\fde)$ be a matched pair of groupoids. There is an associated {\it
diagonal groupoid} $\Do \rightrightarrows \Pc$ with arrows set
$\Do = \Vc _{\,\tgt\!}\times_{\sou} \Hc$, and source, target,
composition and identity given by
\begin{align*} &\sou(g,x)= \sou(g),\quad  \tgt(g,x)= \tgt(x),\\
&(g,x)(h,y)=(g(x\fde h),(x\gde h)y), \quad \mathbf{\id}_P
=(\id_P,\id_P),
\end{align*}
$g,h\in \Vc$, $x,y\in \Hc$, $P\in \Pc$. The diagonal groupoid
$\Do$ admits an \emph{exact factorization} as the product of its
wide subgroupoids $\Vc$ and $\Hc$. In what follows we shall use
the notation $\Do = \Vc\bowtie \Hc$.

 \subsection{Characterization} Out of each matched pair
of groupoids $(\Hc, \Vc)$ one can build a vacant double groupoid
 with  boxes $\B := \Hc {\,}_r\times_t \Vc$. The element
$X = (x,g) \in \Hc {\,}_r\times_t \Vc$ can be represented by $X =
\begin{matrix} \quad \quad x \quad
\\  x \fde g \, \boxee \,\, g \\ \quad \quad x \gde g  \quad
\end{matrix}$.

More precisely, the following result due to Mackenzie says that
the three notions discussed in this section are indeed equivalent
ones.

\begin{prop}\label{equiv-matchedpair} \cite[Theorems 2.10 and 2.15]{mk1}
The following notions are equivalent.

\begin{enumerate}
\item Matched pairs of groupoids. \item Groupoids with an exact
factorization. \item Vacant double groupoids.
\end{enumerate}\qed
\end{prop}

We refer the reader to \cite{AN1} for more details on these
equivalences.

\section{Weak Hopf algebras arising from vacant double
groupoids}\label{weak-hopf}

Throughout this section we shall consider a finite \emph{vacant}
double groupoid $\Tc$. We recall the weak Hopf algebra introduced
in \cite{AN1} arising from $\Tc$.

\medbreak Let us first recall the definition of $\Opext(\Vc,
\Hc)$. As a set $\Opext(\Vc, \Hc)$ consists of equivalence classes
of pairs $(\sigma, \tau)$, where $\sigma$ is a normalized
2-cocycle for the vertical groupoid ${\mathcal B}\rightrightarrows
\Vc$ and $\tau$ a normalized 2-cocycle for the horizontal groupoid
${\mathcal B}\rightrightarrows \Hc$, such that
\begin{equation}
\label{cociclo-sigma-tau}   \quad \sigma(AB, CD) \tau
\left(\begin{matrix} A \\ C \end{matrix},
\begin{matrix}B \\ D \end{matrix}\right)
= \tau(A, B) \tau(C, D) \sigma(A, C) \sigma(B, D),  \end{equation}
for all appropriately composable $A, B, C, D \in {\mathcal B}$.

\medbreak The set $\Opext(\Vc, \Hc)$ has a group structure and it
fits into the following generalization of the so called \emph{Kac
exact sequence}, see \cite{AN1}, \cite{AM}:
\begin{equation}\label{kac}
\begin{aligned}
0 &\longrightarrow  H^1(\D, {\ku^{\times}}) \stackrel{\text{
res}}{\longrightarrow} H^1(\Hc, {\ku^{\times}}) \oplus H^1(\Vc,
{\ku^{\times}})
\longrightarrow \Aut(\ku\,\Tc) \\
&\longrightarrow  H^2(\D, {\ku^{\times}}) \stackrel{\text{
res}}{\longrightarrow} H^2(\Hc, {\ku^{\times}}) \oplus
H^2(\Vc, {\ku^{\times}})\longrightarrow \Opext(\Vc,\Hc) \\
&{\longrightarrow} H^3(\D, {\ku^{\times}})\stackrel{\text{ res
}}{\longrightarrow} H^3(\Hc,{\ku^{\times}}) \oplus H^3(\Vc,
{\ku^{\times}})\longrightarrow \dots,
\end{aligned}
\end{equation}

\noindent where $res$ denotes the restriction maps.

\medbreak
\begin{rmk} Let $(\sigma, \tau)$ be a representative of a class in $\Opext(\ku\Tc)$.
The map $\Opext(\ku\Tc)\to H^3(\Do, \ku^{\times})$ in the above
sequence takes the class of  $(\sigma, \tau)$ to the class of the
3-cocycle $\omega = \omega(\sigma, \tau)$ defined by
\begin{equation}\label{3co} \omega((g, x),(h, y),(f, z))=
\tau(\begin{matrix} x \gde h \quad
\\\boxee \, \, y\fde f \\ \quad \quad \quad
\end{matrix}, \begin{matrix} \quad y \quad
\\ \boxee \,\, f \\ \quad \end{matrix}) \, \,
\sigma(\begin{matrix} x \quad
\\\boxe \, \, h \\  \quad
\end{matrix}, \begin{matrix} x \gde h  \quad
\\ \, \boxe \,\, y\fde f \\ \quad \end{matrix}),
\end{equation}
for all composable $(g, x), (h, y), (f, z)\in \Vc\bowtie \Hc$.
Compare with \cite{Sch}.

\pf We only sketch the proof. The Kac exact sequence \eqref{kac}
comes from a short exact sequence of double complexes $$0\to
{\mathcal A}^{\bullet, \bullet}(\ku^{\times})\to \Bg^{\bullet,
\bullet}(\ku^{\times})\to \mathcal{E}^{\bullet,
\bullet}(\ku^{\times})\to 0, $$ see \cite{AN1}, \cite{AM}. The
associated total complex of the double complex $\Bg^{\bullet,
\bullet}(\ku^{\times})$ is a ``non-standard" projective resolution
of the trivial $\ku\Do$-module, hence there are isomorphisms
$\zeta_n:H^n(\text{Tot}\Bg(\ku^{\times}))\stackrel{\simeq}{\to
}H^n(\Do,\ku^{\times})$ for all $n\geq 0$, that can be computed
explicitly. The map $\Opext(\Vc,\Hc)\to H^3(\D, \ku^{\times})$ is
obtained by composing the map $\Opext(\Vc,\Hc) =
H^2(\text{Tot}{\mathcal A}^{\bullet,
\bullet}(\ku^{\times}))\rightarrow H^3(\text{Tot}\Bg^{\bullet,
\bullet}(\ku^{\times}))$ associated to the inclusion ${\mathcal
A}^{\bullet, \bullet}(\ku^{\times})\to \Bg^{\bullet,
\bullet}(\ku^{\times})$ with $\zeta_3$.\epf \end{rmk}

We note for future use the following properties of the Kac
3-cocycle $\omega = \omega(\sigma, \tau)$:
\begin{flalign}\label{kac11} &\omega\vert_{\Vc \times \Do \times \D} = 1,&\\
\label{kac12} &\omega((g, x), (h, y), (f, z)) = \omega
((\id_{\sou(x)}, x),(h, y), (f, \id_{\tgt(f)})),
&\\
\label{kac13} & \omega((\id_{\sou(x)}, x), (h, \id_{\tgt(h)}), (f,
\id_{\tgt(f)})) = \sigma(\begin{matrix} \quad x \quad
\\  \, \boxee \, h \,  \\ \quad \quad  \quad
\end{matrix}, \begin{matrix} x\gde h
\\  \quad \boxee \,\, f \\ \quad \quad  \quad
\end{matrix}),&\end{flalign}
for all appropriate choice of $(g, x), (h, y), (f, z)\in\Do.$

\medbreak Formulas \eqref{kac12} and \eqref{kac13} are a
consequence of \eqref{3co}, and \eqref{kac11} follows from Remark
3.7 in \cite{AN1}.

\medbreak We now fix a representative $(\sigma, \tau)$ of a class
in $\Opext(\Vc, \Hc)$. We denote by $\ku^{\sigma}_{\tau}\Tc$ the
vector space spanned by the boxes of $\Tc$ endowed with the
$\sigma$-twisted groupoid algebra of the vertical composition
groupoid and the dual $\tau$-twisted groupoid algebra of the
horizontal composition groupoid.

\medbreak That is, the comultiplication and multiplication in
$\ku^{\sigma}_{\tau}\Tc$ are determined by:
\[\Delta(A) = \sum_{A = BC} \tau (B,
C) \quad B\otimes C; \quad A.B = \begin{cases}\sigma(A, B)
\begin{matrix}A\\ B\end{matrix}, \quad \text{ if } \displaystyle
\frac{A}{B},\\ 0, \quad \text{otherwise}. \end{cases}\]

\begin{teo} (\cite[Theorem 3.8]{AN1}) $\ku^{\sigma}_{\tau}
\Tc$ is a semisimple weak Hopf algebra with antipode given by
$$\Ss (A) = \tau(A, A^h)^{-1} \, \sigma(A^{-1}, A^h)^{-1}\, A^{-1}.$$

The source and target subalgebras are, respectively,  the
subspaces spanned by $(\uno_P)_{P \in \Pc}$ and $({}_P\uno)_{P \in
\Pc}$; so they are commutative of dimension $\vert \Pc \vert$.
\qed \end{teo}

\bigbreak The weak Hopf algebra structure on $\ku^{\sigma}_{\tau}
\Tc$ gives rise to a \emph{multi-fusion} category structure on the
category $\Rep(\ku^{\sigma}_{\tau} \Tc)$ of its finite-dimensional
representations. This category is a \emph{fusion} category if and
only if the groupoid $\Vc \rightrightarrows \Pc$ is connected,
c.f. \cite[Proposition 3.11]{AN1}.

\medbreak The tensor category $\Rep(\ku^{\sigma}_{\tau} \Tc)$ can
be described in terms of the combinatorics of double groupoids, as
in \cite[3.4]{AN1}:  objects are $\Hc$-graded vector spaces
endowed with a left $\sigma$-twisted action of the vertical
composition groupoid $\B \rightrightarrows \Hc$ by linear
isomorphisms.

This means that every $A \in \B$ determines a linear isomorphism
$A: V_{b(A)} \to V_{t(A)}$, compatible with vertical composition
in the following sense:
$$A.(B.v) = \sigma(A, B) \, \begin{matrix}A \\ B\end{matrix}.v, $$
for all vertically composable $A, B \in \B$, and for all $v \in
V_{b(A)}$.

\medbreak Tensor product is as in $\Hc$ graded vector spaces, with
twisted $\Vc$-action induced by $\Delta$. The unit object is the
target subalgebra $(\ku^{\sigma}_{\tau}  \Tc)_t = \oplus_{P \in
\Pc} \ku {}_P\uno$, with $\Hc$-grading defined by
$$
(\ku \Tc_t)_x = \begin{cases} \ku {}_{P}\uno, \quad &\text{if } x
= \id P, \\
0, \quad &\text{if $x \notin \Pc$,}
\end{cases}
$$
for all  $x \in \Hc$, and $\B$-action $A . {}_P\uno = \epsilon_t(A
\, {}_P\uno)$.

 \subsection{ The category $\Rep(\ku^{\sigma}_{\tau}
\Tc)$}

We aim to give a combinatorial description of
$\Rep(\ku^{\sigma}_{\tau} \Tc)$ in terms of the matched pair
$(\Hc, \Vc)$ of groupoids associated to $\Tc$.

\medbreak Define the category $\Vect^{\Hc}_{\Vc}(\sigma,\tau)$ as
follows. Objects in $\Vect^{\Hc}_{\Vc}(\sigma,\tau)$ are
$\Hc$-graded vector spaces $M=\oplus_{x\in\Hc} M_x, $ together
with a right action of $\Vc$ by linear isomorphisms $$g: M_x \to
M_{x \gde g},$$ for all $x \in \Hc$ such that $r(x) = t(g)$.

\medbreak That is,  for all homogeneous elements $m\in M$ and all
composable $g,h\in\Vc$ such that $r(|m|)= t(g)$, and for any
$P\in\Pc$

\begin{align}\label{tw1}
m\leftharpoonup \id_P= \begin{cases} m \quad &\text{if } P = r(|m|)\\
0 \quad &\text{otherwise, }
\end{cases}
\end{align}
\begin{equation}\label{tw2} (m\leftharpoonup g)\leftharpoonup h=
\sigma(\begin{matrix} \quad |m| \quad
\\  \, \boxee \,\, g \\ \quad \quad  \quad
\end{matrix}, \begin{matrix} |m|\gde g
\\  \quad \boxee \,\, h \\ \quad \quad  \quad
\end{matrix})\quad m\leftharpoonup (gh),
\end{equation}
\begin{equation}\label{tw3} |m\leftharpoonup g|=|m|\gde g.
\end{equation}

\medbreak The tensor product of two objects $M=\oplus_{x\in\Hc}
M_x, N=\oplus_{y\in\Hc} N_y$ is given by $$(M\otimes
N)_z=\oplus_{xy=z} M_x\otimes N_y,$$ with right action of $\Vc$ on
homogeneous elements
\begin{align}\label{tw4}(m\otimes n)\leftharpoonup g =
\tau(\begin{matrix} |m| \quad \quad
\\  \, \boxee \, |n| \fde g \\ \quad \quad  \quad
\end{matrix}, \begin{matrix} |n|
\\  \quad \boxee \,\, g \\ \quad \quad  \quad
\end{matrix})
 \quad m \leftharpoonup (|n|\fde g) \otimes \;
n\leftharpoonup g.
\end{align}

The unit object is $\ku\Pc$ with $\Hc$-grading determined by
$|P|=\id_P$, $P \in \Pc$, and right $\Vc$-action given by
$P\leftharpoonup g = b(g),$ for all $P\in \Pc$, $g\in\Vc$ such
that $t(g)=P$.

\medbreak

\begin{prop}\label{equi} For any matched pair of groupoids $(\Hc,\Vc)$ and
any pair $(\sigma,\tau)\in \Opext(\ku\Tc)$ there is a natural
tensor equivalence $$\Vect^{\Hc}_{\Vc}(\sigma,\tau)\cong
\Rep(\ku^{\sigma}_{\tau} \Tc^{op}).$$
\end{prop}

\begin{proof} We will actually prove that $\Vect^{\Hc}_{\Vc}(\sigma,\tau)$
is tensor equivalent to the category of right $\ku^{\sigma}_{\tau}
\Tc$-modules.

\medbreak Let  $M$ be a right $\ku^{\sigma}_{\tau} \Tc$-module.
The $\Hc$-grading on $M$ is defined by
\begin{equation}\label{grad}M_x:=M\cdot \begin{matrix} \quad x \quad \\ \begin{tabular}{||p{0,1cm}||} \hline \\
\hline \end{tabular} \\
\quad x \quad
\end{matrix},\end{equation}
for any $x\in\Hc$. If $x\in \Hc$, $h\in\Vc$ are elements such that
$r(x) = t(h)$, then the action of $h$ on $m\in M_x$ is given by
$$ m\leftharpoonup h:= m \cdot \begin{matrix} \quad x \quad
\\  \,\, \boxe \,\, h  \\ \quad \, \quad
\end{matrix}.$$
Straightforward calculations show that equations \eqref{tw1},
\eqref{tw2}, \eqref{tw3} and \eqref{tw4} are fulfilled.

\medbreak Conversely, assume that $N=\oplus_{y\in\Hc} N_y $ is an
object in $\Vect^{\Hc}_{\Vc}(\sigma,\tau)$. Define
\[n \cdot \begin{matrix} \quad x \quad
\\  \,\, \boxe \,\, g  \\ \quad \, \quad
\end{matrix}:=\begin{cases} n \leftharpoonup g \quad &\text{if } y=x\\
 0 \quad &\text{otherwise,}
\end{cases} \quad \forall n\in N_y, \, \begin{matrix} \qquad x \quad
\\  \,\, \boxe \,\, g  \\ \quad \, \quad
\end{matrix}\in \B. \]
This defines a right action and the action on a tensor product is
given by the comultiplication of $\ku^{\sigma}_{\tau} \Tc$.
Moreover, the above functors are strict tensor functors and define
inverse tensor equivalences of categories. Details are left to the
reader.
\end{proof}

\section{$\Rep(\ku^{\sigma}_{\tau} \Tc)$ as a category of bimodules}\label{main-res}

The main goal of this section is to prove that the category
$\Rep(\ku^{\sigma}_{\tau} \Tc)$ is group-theoretical. To achieve
this we shall use an auxiliary tensor category, namely the
category $\ca(\Do,\omega,\Vc)$ introduced in Section
\ref{primera}.

Along this section $\Do \rightrightarrows \Pc$ will be the
diagonal groupoid $\Do = \Vc\bowtie \Hc$ associated to the matched
pair as in Section \ref{vacantes}. We do not assume however that
the wide subgroupoid $\Vc \rightrightarrows \Pc$
 is connected.

\medbreak

We begin with the following technical lemma. For any $P\in\Pc$ let
$$\theta(P)=\# \{g\in\Vc: \sou(g) = P\} = \# \{g\in\Vc:
\tgt(g) = P\}.$$ Observe that $\theta(P)$ is constant on the
connected components of $\Pc$.

Let us define elements $ \Lambda_P, \widetilde{\Lambda}_P,
\Lambda$ in the groupoid algebra $\ku\Vc$ as follows
\begin{equation}\label{lambda}\Lambda_P=\frac{1}{\theta(P)}\sum_{g\in\Vc: \sou(g)=P} g,
\quad \widetilde{\Lambda}_P=\frac{1}{\theta(P)}\sum_{g\in\Vc:
\tgt(g)=P} g, \quad \Lambda=\sum_{P\in\Pc}
\Lambda_P.\end{equation} Note that also $\Lambda = \sum_{P \in
\Pc} \widetilde{\Lambda}_P$.

\begin{lema}\label{integrals}  The following
identities hold, for all $h\in\Vc$:
\begin{itemize}
    \item[(i)] $h\,\Lambda_P=\begin{cases} \Lambda_{\sou(h)} \quad &\text{if } \tgt(h)=P\\
 0 \quad &\text{ otherwise, }
    \end{cases}\,\,$ and $\,\,\widetilde{\Lambda}_P\, h=\begin{cases}\widetilde{\Lambda}_{\tgt(h)} \quad &\text{if } \sou(h)=P\\
 0 \quad &\text{ otherwise, }
    \end{cases}$
    \item[(ii)] $ h \,\Lambda= \id_{\sou(h)}\, \Lambda,$
    \item[(iii)] $\Lambda\, h=  \Lambda\, \id_{\tgt(h)}$.
\end{itemize} \end{lema}

\begin{proof} Straightforward. \end{proof}

For any left $\ku\Vc$-module $M$ we denote $$^{\Vc}\! M:= \{m\in
M: g\rightharpoonup m=\id_{\sou(g)}\rightharpoonup m\text{ for any
} g\in\Vc\}.$$

\begin{rmk} Note that for $M\in \ca(\Do,\omega,\Vc)$, we have
$^{\Vc}\! M=\Lambda\rightharpoonup M$. Indeed
$\Lambda\rightharpoonup M\subseteq ^{\Vc}\! M$ follows from Lemma
\ref{integrals} (ii), and if $m\in \;^{\Vc}\!M$ then
$$\Lambda\rightharpoonup m= \sum_{P\in\Pc} \frac{1}{\theta(P)}\sum_{\sou(g)=P} g\rightharpoonup m=
\sum_{P\in\Pc} \id_P\rightharpoonup m=m. $$ \end{rmk}

\medbreak We now state the main result of this section.

\begin{teo}\label{tequi} Let
$(\sigma, \tau)\in \Opext(\ku\Tc)$ and let $\omega =
\omega(\sigma, \tau)$ be the 3-cocycle defined by \eqref{3co}.
Then the categories $\Rep(\ku^{\sigma}_{\tau} \Tc)$ and
$\ca(\Do,\omega,\Vc)$ are tensor equivalent.
\end{teo}

\medbreak
\begin{proof} Set $p:\Do\to \Hc$, $p(g,x)=x$ and $\pi:\Do\to \Vc$, $\pi(g,x)=g$.
Define the functors $\Phi: \Vect^{\Hc}_{\Vc}(\sigma,\tau)\to
\ca(\Do,\omega,\Vc)$, and $\Psi: \ca(\Do,\omega,\Vc)\to
\Vect^{\Hc}_{\Vc}(\sigma,\tau)$ as follows:
$\Phi(W):=\ku\Vc\otimes W$, where the tensor product is the tensor
product in $\ca(\Do,\omega)$.

The $\Vc$-actions and $\Do$-grading are determined by
\[g\rightharpoonup (h\otimes w):=gh\otimes w,\,\,\, (h\otimes w)
\leftharpoonup g:= h(|w|\fde g)\otimes w\leftharpoonup g,\]
\[|h\otimes w|:= h |w|, \] for
any homogeneous element $w\in W$ and  $g, h \in \Vc$ appropriately
composable.

\bigbreak Given $M\in \ca(\Do,\omega,\Vc)$ define
$\Psi(M):=^{\Vc}\! M$. The $\Hc$-grading on $\Psi(M)$ is given by
\begin{equation} ||m||:=p(|m|),
\end{equation}
and the right action of $\Vc$ is the action on $M$ as an object in
$\ca(\Do,\omega,\Vc)$. Let us prove that this functors are well
defined. Let $W\in \Vect^{\Hc}_{\Vc}(\sigma,\tau)$. Identities
\eqref{sw3}, \eqref{sw1}, \eqref{sw4}, \eqref{sw6} and \eqref{sw2}
are easily checked. Let $w \in W$, $h, g, f \in \Vc$ appropriately
composable, then
\begin{align*} ((h\ot w)\leftharpoonup
g)\leftharpoonup f&= h |w|\fde g\ot\, (w\leftharpoonup
g)\leftharpoonup f\\
&=h (|w|\fde g)\, (|w|\gde g)\fde f\ot\, (w\leftharpoonup
g)\leftharpoonup f\\
&= \omega(|w|,g,f)\;\, h (|w|\fde gf) \ot\, w\leftharpoonup gf,
\end{align*}
the last identity by \eqref{tw2} and \eqref{kac13}. This proves
identity \eqref{sw5}, thus $\Phi$ is well defined.

\medbreak We claim that the functors $\Phi,  \Psi$ give an
equivalence of categories.

Let $\Lambda$ be the element defined by \eqref{lambda}. Define the
map
$$\gamma_M:M\to \Phi(\Psi(M)), \;\, \gamma_M(m)=\pi(|m|)\ot\, \Lambda
\rightharpoonup m $$ for any homogeneous element $m\in M$.

The map $\gamma_M$ is a $\ku\Vc$-bimodule map: indeed, if
$g\in\Vc$ and $ m\in M$ is homogeneous of degree $(h,x)\in\Do$
then
\begin{align*} \gamma_M(g\rightharpoonup m)&= \pi(g|m|)\ot\, \Lambda
\rightharpoonup (g\rightharpoonup m) = g\pi(|m|)\ot\, \Lambda\,g
\rightharpoonup
 m\\
&=g\pi(|m|)\ot\, \Lambda\, \id_{\tgt(g)} \rightharpoonup
 m=g\rightharpoonup \gamma_M(m).
\end{align*}
The third equality by Lemma \ref{integrals} (iii). On the other
hand
\begin{align*} \gamma_M(m\leftharpoonup g)&= \pi(|m|g) \ot\,
\Lambda\rightharpoonup (m\leftharpoonup g)\\
&=h (|m| \fde g) \ot\, (\Lambda\rightharpoonup m)\leftharpoonup g
=\gamma_M(m) \leftharpoonup g.
\end{align*}
The inverse of $\gamma_M$ is given by $\overline{\gamma}_M(g\ot
\Lambda\rightharpoonup m)=gh^{-1}\rightharpoonup m,$ if $m\in M$
is homogeneous of degree $(h,x)$.

\begin{claim} The map $\overline{\gamma}_M$ is well defined. \end{claim}

\pf First note that if $m\in M$ is an homogeneous element of
degree $(h,x)$ then the homogeneous component of degree $(\id,x)$
of $\Lambda_{\sou(h)}\rightharpoonup m$ is $h^{-1}\rightharpoonup
m$.

\medbreak

Now, let $m, n\in M$ be homogeneous elements of degree $(h,x)$ and
$(h',x')$ respectively, such that $\Lambda\rightharpoonup
m=\Lambda\rightharpoonup n$. Equation \eqref{sw4} implies that
$\Lambda_{\sou(h)}\rightharpoonup
m=\Lambda_{\sou(h')}\rightharpoonup n$. Any homogeneous component
of $\Lambda_{\sou(h)}\rightharpoonup m$, respectively
$\Lambda_{\sou(h')}\rightharpoonup n$, has degree $(gh,x)$, resp.
$(gh',x')$, for some $g\in\Vc$. This implies that $x=x'$. By the
previous observation $h^{-1}\rightharpoonup
m=h'^{-1}\rightharpoonup n$.\epf

\medbreak If $W \in \Vect^{\Hc}_{\Vc}(\sigma,\tau)$ then
$\Psi(\Phi(W))=^{\Vc}(\ku\Vc\otimes W)\cong W$,  the isomorphism
$\phi_W: W\longrightarrow\, ^{\Vc}(\ku\Vc\otimes W)$, being given
by $\phi(w):= \Lambda\,\id_{\tgt(|w|)}\otimes w$,  for any
homogeneous element $w\in W$.

\medbreak Finally, we shall prove that $\Phi$ is a tensor functor.
Let $W, U \in \Vect^{\Hc}_{\Vc}(\sigma,\tau)$ and define
$\xi_{WU}: \Phi(W\otimes U)\to\Phi(W)\,\otb\, \Phi(U)$ in the form
\[ \xi_{WU}(g\otimes (w\otimes u)):=(g\otimes w)\otb
(\tgt(|w|)\otimes u),\] for any appropriate choice of $g\in \Vc$,
and homogeneous $w\in W$, $u\in U$.

\medbreak The natural map $\xi_{WU}$ is an isomorphism, its
inverse $\overline\xi_{WU}: \Phi(W)\,\otb\, \Phi(U) \to
\Phi(W\otimes U)$ being given by
\[ \overline\xi_{WU}((h \otimes z) \otimes (f \otimes y)) = h(|z| \fde f) \otimes z \leftharpoonup f \otimes y. \]

We next show that $\xi_{WU}$ is a morphism of $\Vc$-bimodules.

Let $h, g \in \Vc$, $w \in W$, $U \in U$. To prove right
$\Vc$-linearity we compute
\begin{align*}&\xi_{WU}((g\otimes (w\otimes
u))\leftharpoonup h)= \xi_{WU}(g (|w||u|\fde h)\otimes (w\otimes
u)\leftharpoonup h)\\
& = \tau(\begin{tiny}\begin{matrix} |w|\quad \quad
\\\boxe \, \, |u|\fde h \\ \quad \quad \quad
\end{matrix}, \begin{matrix}  |u| \quad
\\ \boxe \,\, h \\ \quad \end{matrix}\end{tiny})
\,\, \xi_{WU}(g (|w||u|\fde h)\otimes (w\leftharpoonup (|u|\fde
h)\otimes
u\leftharpoonup h))\\
& = \tau(\begin{tiny}\begin{matrix} |w|\quad \quad
\\\boxe \, \, |u|\fde h \\ \quad \quad \quad
\end{matrix}, \begin{matrix}  |u| \quad
\\ \boxe \,\, h \\ \quad \end{matrix}\end{tiny}) \,\, (g (|w||u|\fde h)\otimes
w\leftharpoonup (|u|\fde h))\otb\, (\sou(|u|) \otimes
u\leftharpoonup h)\\
& = \tau(\begin{tiny}\begin{matrix} |w|\quad \quad
\\\boxe \, \, |u|\fde h \\ \quad \quad \quad
\end{matrix}, \begin{matrix}  |u| \quad
\\ \boxe \,\, h \\ \quad \end{matrix}\end{tiny}) \,\, ((g\otimes w)\leftharpoonup
(|u|\fde h)) \otb\, (\sou(|u|) \otimes u\leftharpoonup h)\\
& = \tau(\begin{tiny}\begin{matrix} |w|\quad \quad
\\\boxe \, \, |u|\fde h \\ \quad \quad \quad
\end{matrix}, \begin{matrix}  |u| \quad
\\ \boxe \,\, h \\ \quad \end{matrix}\end{tiny}) \,\,(g\otimes w)
\otb\, (|u|\fde h\otimes u\leftharpoonup h)\\
&= \omega (g|w|,|u|,h)\,\, (g\otimes w)
\otb\, (|u|\fde h\otimes u\leftharpoonup h)\\
&= \xi_{WU}(g\otimes (w\otimes u))\leftharpoonup h.
\end{align*}
The fifth equality by \eqref{tp1}, the sixth follows from the
definition of $\omega$, and the last equality follows from
equation \eqref{ap-asoc}. Left $\Vc$-linearity is similarly
established.

\begin{claim} \label{tensor} For all $U,V,W
\in \Vect^{\Hc}_{\Vc}(\sigma,\tau)$, we have
$$a_{\Phi(U),\Phi(V),\Phi(W)} \, (\xi_{U,V}\ot \id) \, \xi_{U\ot V, W} =
(\id\ot\xi_{V,W}) \, \xi_{U,V\ot W} \, \Phi(a_{U,V,W}).$$
\end{claim}

\pf Let $U,V,W \in \Vect^{\Hc}_{\Vc}(\sigma,\tau)$,  $g\in \Vc$
and let $u\in U, v\in V, w\in W$   be homogeneous elements with
appropriately composable degree. The left hand side of this
equation evaluated in $g\ot u\ot v\ot w$ gives
\begin{align*} a\,(\xi_{U,V}\ot \id)&(g\ot
(u\ot v)\otb \, \sou(|w|) \ot w) = a((g\ot u\otb\,
\sou(|v|)\ot v)\otb \, \sou(|w|)\ot w)\\
& = \omega(g|u|,|v|,|w|)\,\, g\ot u\otb\,( \sou(|v|) \ot v\otb \,
\sou(|w|) \ot w)\\ & = g\ot u\otb\,(\sou(|v|)\ot v\otb \,
\sou(|w|)\ot w).
\end{align*}
The last equality by \eqref{kac12}. The right hand side of
\eqref{tensor} evaluated in $g\ot u\ot v\ot w$ gives
\begin{align*} (\id\ot\xi_{V,W}) \xi_{U,V\ot W} (g\ot u\ot (v\ot
w))&=(\id\ot\xi_{V,W})(g\ot u\otb\, \sou(|v|)\ot (v\ot w))\\
&=g\ot u \otb\,(\sou(|v|)\ot v \otb\, \sou(|w|)\ot w).
\end{align*} \epf
The proof of the theorem is now complete. \end{proof}

We now can state the main result of this section.

\begin{teo}\label{gp-ttic} Suppose that $\Rep(\ku^{\sigma}_{\tau} \Tc)$ is a fusion
category. Let $\Do$ be the diagonal groupoid associated to $\Tc$,
and let $\omega = \omega(\sigma, \tau) \in H^3(\Do, \ku^{\times})$
be the 3-cocycle given by \eqref{3co}. Then there is an
equivalence of tensor categories
$$\Rep(\ku^{\sigma}_{\tau} \Tc) \simeq
\ca(D, \overline{\omega}, V),$$ where $\overline \omega$ is a
normalized 3-cocycle on $D$ cohomologous to the restriction
$\widehat \omega$.

In particular, the category $\Rep(\ku^{\sigma}_{\tau} \Tc)$ is
group-theoretical. \qed
\end{teo}

\pf By \cite[Proposition 3.11]{AN1}, the assumption implies that
$\Vc \rightrightarrows \Pc$ is connected. Combining Corollary
\ref{equiv-bimod} and Theorem \ref{tequi}, we get tensor
equivalences
$$\Rep(\ku^{\sigma}_{\tau} \Tc) \simeq
\ca(\Do, \omega, \Vc) \simeq \ca(D, \widehat \omega, V, \widehat
\psi),$$ where $\psi \in H^2(\Vc, \ku^{\times})$ is determined by
\eqref{coborde}, and $\widehat \psi$ is the 2-cocycle on $V$
obtained by restriction. By \cite[Remark 3.2]{fs-indic}, there
exists a 3-cocycle $\overline \omega$ on $D$ which is cohomologous
to $\widehat \omega$ and such that $\ca(D, \widehat \omega, V,
\widehat \psi) \simeq \ca(D, \overline \omega, V)$. This finishes
the proof of the theorem. \epf

It follows that $\Rep(\ku^{\sigma}_{\tau} \Tc)$ is the
representation category of certain (unique up to gauge
equivalence) semisimple finite-dimensional quasi-Hopf algebra. We
point out that an explicit description (up to twist equivalence)
of this quasi-Hopf algebra has been given in \cite{fs-indic},
where the Frobenius-Schur indicators of a group-theoretical fusion
category were also computed.

By specializing the results of Ostrik on classification of module
categories over group-theoretical categories, we can parameterize
module categories and fiber functors for the categories
$\Rep(\ku^{\sigma}_{\tau} \Tc)$.

\bibliographystyle{amsalpha}

\end{document}